\documentclass[]{preprint}

\usepackage{subfigure}

\theoremstyle{plain}
\newtheorem{theorem}{Theorem}[section]

\newtheorem{definition}[theorem]{Definition}

\theoremstyle{remark}

\begin{document}




\title{{\itshape Preprocessing of centred logratio transformed density functions using smoothing splines}}

\author{J. Machalov\'a$^{\rm a}$$^{\ast}$,
\thanks{$^\ast$Corresponding author. Email: jitka.machalova@upol.cz
\vspace{6pt}} K. Hron$^{\rm ab}$,
\vspace{6pt} G.S. Monti$^{\rm c}$\\
\vspace{6pt}  $^{\rm a}${\em Department of Mathematical Analysis and Applications of Mathematics, Faculty of Science, Palack\'y University, 17. listopadu 12, CZ-77146 Olomouc, Czech Republic,}\\
$^{\rm b}${\em Department of Geoinformatics, Faculty of Science, Palack\'y University, t\v r. Svobody 26, CZ-77146 Olomouc, Czech Republic,}\\
$^{\rm c}${\em Department of Economics, Management and Statistics, University of Milano-Bicocca, Via Bicocca degli Arcimboldi 8,
20126 Milano, Italy}\\
}

\maketitle

\begin{abstract}
With large-scale database systems, statistical analysis of data, formed by probability distributions, become an important task in explorative data analysis. Nevertheless, due to specific properties of density functions, their proper statistical treatment still represents a challenging task in functional data analysis. Namely, the usual $L^2$ metric does not fully accounts for the relative character of information, carried by density functions; instead, their geometrical features are followed by Bayes spaces of measures. The easiest possibility of expressing density functions in $L^2$ space is to use centred logratio transformation, nevertheless, it results in functional data with a constant integral constraint that needs to be taken into account for further analysis. While theoretical background for reasonable analysis of density functions is already provided comprehensively by Bayes spaces themselves, preprocessing issues still need to be developed. The aim of this paper is to introduce optimal smoothing splines for centred logratio transformed density functions that take all their specific features into account and provide a concise methodology for reasonable preprocessing of raw (discretized) distributional observations. Theoretical developments are illustrated with a real-world data set from official statistics.
\end{abstract}


\begin{classcode}
\textit{Classification codes}: 62H99, 65D07, 65D10
\end{classcode}


\section{Introduction}


Density functions, i.e. Borel measurable, positive functions on a support $I$ with a unit integral constraint, can be considered as a special case of functional data \cite{ramsay05}. Nowadays, they frequently occur in practice due to large-scale database systems, where the individual observations are summarized by distributions to preserve their intrinsic variability and enable to analyze statistically groups of individuals in meaningful way \cite{billard07, noirhomme11}, i.e. to form distribution-valued variables. Nevertheless, for density functions the usual $L^2$ metric seems to be not appropriate. Not just because of the unit-integral constraint (that forms rather a proper representation due to probability conventions than an inherent feature of densities themselves), in addition, density functions also contain information on relative contributions of Borel sets of real line to the overall probability on the (possibly infinite) support of the corresponding random variables. As a way out, Bayes spaces with se\-pa\-rable Hilbert space properties were proposed for geometric representation of density functions \cite{egozcue06, boogaart10, boogaart14}. In order to enable standard statistical analysis of density functions, an isomorphic mapping from the Bayes space to the standard $L^2$ space is needed. Considering the finite interval support case (specially $I=[a,b]$ for real $a<b$), where Lebesgue measure as reference measure can be used \cite{egozcue06}, the simplest such isomorphism is represented by the \textit{centred logratio (clr) transformation}, defined for a density function $f(x)$ as
\begin{equation}\label{fclr}
    \mbox{clr}[f(x)]=f_c(x)=\ln f(x)-\frac{1}{\eta}\int_I\ln f(x)\,\mbox{d}x,
\end{equation}
where $\eta$ stands for length of the interval $I$, in particular $\eta=b-a$. Because clr transformation forms a one-to-one mapping, the inverse clr transformation is obtained as
\begin{equation}\label{fclrinv}
    \mbox{clr}^{-1}[f_c(x)]=\frac{\exp(f_c(x))}{\int_I\exp(f_c(x))\,\mbox{d}x};
\end{equation}
the denominator is used just to achieve the unit integral constraint representation of the resulting density (without loss of relative information, carried by the density function). Equations (\ref{fclr}) and (\ref{fclrinv}) come as a functional version of the popular centred logratio transformation of compositional data, multivariate observations that carry relative information, see \cite{aitchison86, pawlowsky11} for details. Nevertheless, because obviously
\begin{equation}\label{fcond}
    \int_I\mbox{clr}[f(x)]\,\mbox{d}x=0,
\end{equation}
this additional condition needs to be taken into account for computation and analysis of clr transformed density functions. Although the clr transformation may lead to computational problems for some of statistical methods due to constraint (\ref{fcond}), it is still well acceptable for distance-based methods or functional principal component analysis, similarly as for the case of compositional data, and could thus extend the currently existing analytical tools \cite{delicado11, menafoglio14}.

However, density functions (as well as functional data in general) occur in the practice rarely in their continuous form. For example, the aggregation of individual observations in case of the distribution-valued variables leads naturally to discretized form of histogram data that need to be approximated (smoothed) by an appropriate function. No wonder that approximation of nonparametric distributions is one of basic problems of functional data analysis \cite{ramsay05} and a number of publications are devoted to cope with this issue \cite{bowman97, silverman86}. For the purpose of functional data analysis, smoothing $B$-splines turned out to be the most appropriate approximative tool as the resulting coefficients of basis functions can be directly used for further statistical analysis \cite{ramsay05}. Nevertheless, taking the inherent features of densities into account, the possibility of obtaining smoothing splines (and even $B$-splines) in case of density functions gets quite a complex problem \cite{dias98, eilers96, gu93, maechler95}. In \cite{maechler96} logarithmic transformation is proposed to simplify the estimation process and properties of the resulting splines, nevertheless, without any deeper methodological background.

Although the methodology of Bayes spaces was successfully applied to theoretical problems related with Bayesian approach to statistical analysis \cite{egozcue13,boogaart10}, its application to statistical processing of density functions is still limited due to absence of a reasonable approximation tool that would enable to proceed from functional data to smooth functions. Neither smoothing the original discretised densities \cite{delicado11} nor using of Bernstein polynomials, that is proposed in \cite{menafoglio14}, is coherent with the Bayes space methodology.

A natural and logical step to perform smoothing of density functions using $B$-splines is to express the data in the $L^2$ space of clr transformed densities and to perform the computations there. Consequently, it is guaranteed that all theoretical properties of smoothing $B$-splines hold and no additional reasoning is necessary. Nevertheless, the condition (\ref{fcond}) needs to be taken into account for estimation of spline coefficients. The aim of the paper is to provide a concise solution of this problem from the point of view of numerical mathematics.

The paper is organized as follows. The next section is devoted to general theory of $B$-spline representation, followed by the optimal smoothing problem, representing a trade-off between interpolation and the least squares approximation, in Section 3. In this section also a solution for the conditional smoothing is proposed. Section 4 continues with concrete example, based on a real-world data set; finally, Section 5 concludes.


\section{The $B$-spline representation}


In this section we recall some basics concerning the $B$-spline representation of splines, see for example \cite{deboor78,dierckx93}. Let the sequence of knots $\Delta\lambda:=\lambda_{0}=a<\lambda_{1}<\ldots<\lambda_{g}<b=\lambda_{g+1}$  be given. In the following, ${\cal S}_{k}^{\Delta\lambda}[a,b]$ denotes the vector space of polynomial splines of degree $k>0$, defined on a finite interval $[a,b]$ with the sequence of knots $\Delta\lambda$. It is known that $dim\left({\cal S}_{k}^{\Delta\lambda}[a,b]\right)=g+k+1$. Then every spline $s_{k}(x)\in{\cal S}_{k}^{\Delta\lambda}[a,b]$ has a unique representation
\begin{equation}\label{br}
    s_{k}\left(x\right)=\sum_{i=-k}^{g}b_{i}B_{i}^{k+1}\left(x\right).
\end{equation}
Vector ${\bf b}=(b_{-k},\ldots,b_{g})^{\top}$ is {\it vector of $B$-spline coefficients} of~$s_{k}(x)$, functions $B_{i}^{k+1}\left(x\right)$, $i=-k,\ldots,g$ are {\it $B$-splines of degree $k$} and form basis in ${\cal S}_{k}^{\Delta\lambda}[a,b]$.

\begin{definition}
Let ${\bf x}=(x_{1},\ldots,x_{n})^{\top}$ be given and let $\bigl\{B_{i}^{k+1}(x)\bigr\}_{i=-k}^{g}$ be {\it B}-spline basis of~${\cal S}_{k}^{\Delta\lambda}[a,b]$. Then
\begin{equation}\label{kol}
    {\bf C}_{k+1}({\bf x})=\left(\begin{array}{ccc}
                                   B_{-k}^{k+1}\left(x_{1}\right) & \ldots   & B_{g}^{k+1}\left(x_{1}\right)\\
                                                           \vdots & \ddots & \vdots\\
                                   B_{-k}^{k+1}\left(x_{n}\right) & \ldots   & B_{g}^{k+1}\left(x_{n}\right)
                                 \end{array}\right)\in\mathbb{R}^{n,g+k+1}
\end{equation} is called {\it the collocation matrix}.
\end{definition}

Consequently, every spline from ${\cal S}_{k}^{\Delta\lambda}[a,b]$ can be written in matrix notation as
\begin{equation}\label{mr}
    s_{k}(x)={\bf C}_{k+1}(x){\bf b}.
\end{equation}

Now let $l\in\left\{1,\ldots,k-1\right\}$. It is known that derivative of order $l$ of the spline $s_{k}(x)\in{\cal S}_{k}^{\Delta\lambda}[a,b]$ is a spline $s_{k-l}(x)\in{\cal S}_{k-l}^{\Delta\lambda}[a,b]$ with the same knots. Using properties of $B$-splines the spline derivatives can be written in matrix notation as
$$s_{k}^{(l)}\left(x\right)={\bf C}_{k+1-l}(x){\bf b}^{(l)},$$
where ${\bf b}^{(l)}\in\mathbb{R}^{g+k+1-l}$ is given by
\begin{eqnarray*}
    {\bf b}^{(l)} & = & {\bf D}_{l}{\bf L}_{l}{\bf b}^{(l-1)} \\
                  & = & {\bf D}_{l}{\bf L}_{l}\ldots{\bf D}_{1}{\bf L}_{1}{\bf b}\\
                  & = & {\bf S}_{l}{\bf b}
\end{eqnarray*}
and ${\bf b}^{(0)}={\bf b}$. Upper triangular matrix ${\bf S}_{l}={\bf D}_{l}{\bf L}_{l}\ldots{\bf D}_{1}{\bf L}_{1}\in\mathbb{R}^{g+k+1-l,g+k+1}$ has full row rank;
hereat ${\bf D}_{j}\in\mathbb{R}^{g+k+1-j,g+k+1-j}$ is diagonal matrix such that
$${\bf D}_{j}=\left(k+1-j\right)diag\left(d_{-k+j},\ldots,d_{g}\right)$$
with $$d_{i}=\dfrac{1}{\lambda_{i+k+1-j}-\lambda_{i}}\quad\forall i=-k+j,\ldots,g$$
and $${\bf L}_{j}:=\left(\begin{array}{cccc}
                           -1 & 1      &        &  \\
                              & \ddots & \ddots &  \\
                              &        & -1     & 1
                         \end{array} \right)\in\mathbb{R}^{g+k+1-j,g+k+2-j}.$$


\section{The optimal smoothing problem}


The optimal smoothing problem represents a compromise between the interpolation problem and the least squares approximation. Let data $(x_{i},y_{i})$, $a\leq x_{i}\leq b$, weights $w_{i}\geq 0$, $i=1,\ldots,n$, $n\geq g+1$ and parameter $\alpha>0$ are given. For $l\in\left\{1,\ldots,k-1\right\}$ the task is to find a spline $s_{k}(x)\in{\cal S}_{k}^{\Delta\lambda}[a,b]$, which minimizes functional
\begin{equation}\label{jl}
    J_{l}(s_{k})=\int_{a}^{b}\left[s_{k}^{(l)}(x)\right]^{2}\,\mbox{d}x+\alpha\sum\limits_{i=1}^{n}w_{i}\left[y_{i}-s_{k}(x_{i})\right]^{2}.
\end{equation}
This spline is called {\it the smoothing spline.}
Further we denote $\mathbf{x}=\left(x_{1},\ldots,x_{n}\right)^{\top}$, $\mathbf{y}=\left(y_{1},\ldots,y_{n}\right)^{\top}$, $\mathbf{w}=\left(w_{1},\ldots,w_{n}\right)^{\top}$ and $\mathbf{W}=diag\left(\mathbf{w}\right)$. The functional $J_{l}(s_{k})$ can be written in a matrix form as
\begin{equation}\label{jb}
    J_{l}({\bf b})={\bf b}^{\top}{\bf N}_{kl}{\bf b}+\alpha\left[{\bf y}-{\bf C}_{k+1}({\bf x}){\bf b}]^{\top}{\bf W}[{\bf y}-{\bf C}_{k+1}({\bf x}){\bf b}\right].
\end{equation}
Matrix ${\bf N}_{kl}={\bf S}_{l}^{\top}{\bf M}_{kl}{\bf S}_{l}$ is positive semidefinite, where
$${\bf M}_{kl}=\left(\begin{array}{ccc}
                       \left(B_{-k+l}^{k+1-l},B_{-k+l}^{k+1-l}\right) & \ldots  & \left(B_{g}^{k+1-l},B_{-k+l}^{k+1-l}\right)\\
                                                               \vdots &         & \vdots\\
                          \left(B_{-k+l}^{k+1-l},B_{g}^{k+1-l}\right) & \ldots  & \left(B_{g}^{k+1-l},B_{g}^{k+1-l}\right)
                     \end{array}\right)\in\mathbb{R}^{g+k+1-l,g+k+1-l}$$
and
$$\left(B_{i}^{k+1-l},B_{j}^{k+1-l}\right)=\int\limits_{a}^{b}B_{i}^{k+1-l}(x)B_{j}^{k+1-l}(x)\,\mbox{d}x$$
stands for scalar product of $B$-splines in $L^2([a,b])$ space. Matrix ${\bf M}_{kl}$ is positive definite, because $B_{i}^{k+1-l}(x)\geq 0$, $i=-k+l,\ldots,g$ are basis functions.

\bigskip
Our aim is to find a spline $s_{k}(x)\in{\cal S}_{k}^{\Delta\lambda}[a,b]$, which minimizes functional $J_{l}(s_{k})$, in other words, we want to find a minimum of function $J_{l}(\mathbf{b})$. From the necessary and sufficient condition for a minimum of this function, i.e.
$$\frac{\partial J_{l}({\bf b})}{\partial{\bf b}^{\top}}=0,\quad l\in\left\{1,\dots,k-1\right\},$$
we get a system of linear equations
\begin{equation}\label{sm}
    \left[\alpha^{-1}{\bf N}_{kl}+{\bf C}_{k+1}^{\top}({\bf x}){\bf W}{\bf C}_{k+1}({\bf x})\right]{\bf b}={\bf C}_{k+1}^{\top}({\bf x}){\bf Wy}.
\end{equation}
Solution of this system exists if and only if this system is consistent, i.e.
\begin{equation*}
{\bf C}_{k+1}^{\top}({\bf x}){\bf Wy}\in{\cal R}\left(\alpha^{-1}{\bf N}_{kl}+{\bf C}_{k+1}^{\top}({\bf x}){\bf W}{\bf C}_{k+1}({\bf x})\right),\end{equation*}
and it is given by
\begin{equation}\label{b**}
    {\bf b}^{*}=\left[\alpha^{-1}{\bf N}_{kl}+{\bf C}_{k+1}^{\top}({\bf x}){\bf WC}_{k+1}({\bf x})\right]^{-}{\bf C}_{k+1}^{\top}({\bf x}){\bf Wy}.
\end{equation}
Here $\mathcal{R}(\mathbf{A})$ denotes the range of matrix $\mathbf{A}$. Matrix $\mathbf{A}^{-}$ is a generalized inverse of $\mathbf{A}$, i.e. a matrix such that $\mathbf{A}\mathbf{A}^{-}\mathbf{A}=\mathbf{A}$, see \cite{rao71} for details. If matrix $\mathbf{A}$ is regular, then $\mathbf{A}^{-}=\mathbf{A}^{-1}$.

Accordingly, if the matrix of the system (\ref{sm}) is regular, then there exists a unique solution of (\ref{sm}), i.e. there exists a unique smoothing spline.  If matrix of (\ref{sm}) is singular and the system (\ref{sm}) is consistent, then there exists a class of solutions. In this case we can find {\it a minimum norm solution} of (\ref{sm}), i.e. a solution which has the smallest norm among all solutions. Such solution is unique and is given by
\begin{equation}\label{bm**}
    {\bf b}^{*}=\left[\alpha^{-1}{\bf N}_{kl}+{\bf C}_{k+1}^{\top}({\bf x}){\bf WC}_{k+1}({\bf x})\right]^{-}_{m}{\bf C}_{k+1}^{\top}({\bf x}){\bf Wy};
\end{equation}
then we refer to a unique {\it optimal smoothing spline}.
The symbol $\mathbf{A}^{-}_{m}$ denotes {\it a minimum norm generalized inverse}, i.e. a matrix such that
$$\mathbf{A}\mathbf{A}^{-}_{m}\mathbf{A}=\mathbf{A},\qquad\left(\mathbf{A}^{-}_{m}\mathbf{A}\right)^{\top}=\mathbf{A}^{-}_{m}\mathbf{A}.$$

Consequently, the required smoothing spline or optimal smoothing spline are obtained by formula
\begin{equation}\label{sb*}
    s_{k}^{*}(x)=\sum\limits_{i=-k}^{g}b_{i}^{*}B_{i}^{k+1}(x),
\end{equation}
where the vector $\mathbf{b}^{*}=\left(b_{-k}^{*},\ldots,b_{g}^{*}\right)^{\top}$ is given in (\ref{b**}) and (\ref{bm**}), respectively.

\bigskip
In the case of smoothed clr transformed density functions an additional condition
\begin{equation}\label{p}
    \int_{a}^{b}s_{k}(x)\,\mbox{d}x=0
\end{equation}
needs to be fulfilled. From $B$-spline representation of splines it is known that the spline
$$s_{k}(x)=\sum\limits_{i=-k}^{g}b_{i}B_{i}^{k+1}(x)$$
is a derivative of spline
$$s_{k+1}(x)=\sum\limits_{i=-k-1}^{g}c_{i}B_{i}^{k+2}(x),$$
if
\begin{equation}\label{bc}
    b_{i}=\left(k+1\right)\dfrac{c_{i}-c_{i-1}}{\lambda_{i+k+1}-\lambda_{i}}\quad\forall i=-k,\ldots,g.
\end{equation}

\medskip
Now let's consider a spline $s_{k}(x)\in{\cal S}_{k}^{\Delta\lambda}[a,b]$ such that condition (\ref{p}) holds, i.e.
\begin{equation*}
    0=\int_{a}^{b}s_{k}(x)\,\mbox{d}x=\left[s_{k+1}\left(x\right)\right]_{a}^{b}=s_{k+1}(\lambda_{g+1})-s_{k+1}(\lambda_{0}),
\end{equation*}
because $a=\lambda_{0}$, $b=\lambda_{g+1}$. With respect to the definition, properties of $B$-splines and the mentioned additional knots we get
\begin{equation*}
    0=s_{k+1}(\lambda_{g+1})-s_{k+1}(\lambda_{0})=c_{g}-c_{-k-1},
\end{equation*}
so that
\begin{equation}\label{cc}
    c_{-k-1}=c_{g}.
\end{equation}
Now we find a relationship between the vector ${\bf b}=(b_{-k},\ldots,b_{g})^{\top}$ of $B$-spline coefficients of $s_{k}(x)$ and the vector ${\bf c}=(c_{-k-1},\ldots,c_{g})^{\top}$ of $s_{k+1}(x)$, ${\bf c}\in\mathbb{R}^{g+k+2}$.  Using (\ref{cc}) we can write
\begin{equation}\label{bc1}
    \mathbf{b}=\mathbf{D}\mathbf{K}\mathbf{\bar{c}},
\end{equation}
where ${\bf \bar{c}}=(c_{-k},\ldots,c_{g})^{\top}\in\mathbb{R}^{g+k+1}$,
$$\mathbf{D}=\left(k+1\right)diag\left(\dfrac{1}{\lambda_{1}-\lambda_{-k}},\ldots,\dfrac{1}{\lambda_{g+k+1}-\lambda_{g}}\right)\in\mathbb{R}^{g+k+1,g+k+1}$$
and
$$\mathbf{K}=\left(\begin{array}{rrrrr}
	                 1 &  0 & 0 & \cdots & -1 \\
                    -1 &  1 & 0 & \cdots &  0 \\
                     0 & -1 & 1 & \cdots &  0 \\
                     \vdots &  \vdots & \ddots & \ddots & \vdots \\
                     0 & 0  & \cdots & -1 & 1 \\
                   \end{array}\right)\in\mathbb{R}^{g+k+1,g+k+1}.$$
Note that both matrices $\mathbf{D}$ and $\mathbf{K}$ are regular. The goal is to find a vector $\mathbf{b}$, which minimizes function $J_{l}(\mathbf{b})$ given in (\ref{jb}) and satisfies condition (\ref{bc1}) simultaneously. Using both relationships we can rewrite function $J_{l}(\mathbf{b})$ as
\begin{equation*}\label{jc}
    J_{l}({\bf\bar{c}})={\bf\bar{c}}^{\top}\mathbf{K}^{\top}\mathbf{D}^{\top}{\bf N}_{kl}\mathbf{DK\bar{c}}+\alpha\left[{\bf y}-{\bf C}_{k+1}({\bf x})\mathbf{D}\mathbf{K}{\bf\bar{c}}\right]^{\top}{\bf W}\left[{\bf y}-{\bf C}_{k+1}({\bf x})\mathbf{DK\bar{c}}\right].
\end{equation*}
Minimum of $J_{l}(\mathbf{\bar{c}})$ is a solution of the following system,
\begin{equation}\label{sc}
    \left[\alpha^{-1}\left(\mathbf{DK}\right)^{\top}{\bf N}_{kl}\mathbf{DK}+\left({\bf C}_{k+1}({\bf x})\mathbf{DK}\right)^{\top}{\bf W}{\bf C}_{k+1}({\bf x})\mathbf{DK}\right]{\bf
    \bar{c}}=\left({\bf C}_{k+1}({\bf x})\mathbf{DK}\right)^{\top}{\bf Wy}.
\end{equation}

For this system of linear equations we can perform similar considerations as (\ref{sm}). Concretely, if the system (\ref{sc}) is consistent and its matrix is not regular, then a solution of (\ref{sc}) exists,
\begin{equation*}
    {\bf\bar{c}}^{*}=\left[\alpha^{-1}\left(\mathbf{DK}\right)^{\top}{\bf N}_{kl}\mathbf{DK}+\left({\bf C}_{k+1}({\bf x})\mathbf{DK}\right)^{\top}{\bf W}{\bf C}_{k+1}({\bf x})\mathbf{DK}\right]^{-}\mathbf{K}^{\top}\mathbf{D}^{\top}{\bf C}_{k+1}^{\top}({\bf x}){\bf Wy}.
\end{equation*}
The minimum norm solution among all solutions of (\ref{sc}) is given by
\begin{equation*}
    {\bf\bar{c}}^{*}=\left[\alpha^{-1}\left(\mathbf{DK}\right)^{\top}{\bf N}_{kl}\mathbf{DK}+\left({\bf C}_{k+1}({\bf x})\mathbf{DK}\right)^{\top}{\bf W}{\bf C}_{k+1}({\bf x})\mathbf{DK}\right]^{-}_{m}\mathbf{K}^{\top}\mathbf{D}^{\top}{\bf C}_{k+1}^{\top}({\bf x}){\bf Wy}.
\end{equation*}
Finally, the vector of $B$-spline coefficients is obtained as
\begin{equation*}
    \mathbf{b}^{*}=\mathbf{DK}\mathbf{\bar{c}}^{*}.
\end{equation*}

It is important to note that the resulting $B$-spline coefficients fulfill the condition of zero sum, known from discrete version of the clr transformation. This is important for further theoretical reasoning, when functional counterparts to multivariate statistical methods for compositional data (like principal component analysis \cite{delicado11, filzmoser09}) are introduced.

The resulting splines can be used either to back-transform them to the original Bayes space of density functions (in order to provide their smoothed counterparts for visualization and interpretation purposes), or for their further statistical analysis in the clr space. In the next section the smoothing step is performed for real data from an Italian household survey.


\section{Application to real-world data}


In the following, the above theoretical considerations are applied to a real-world data set from official statistics. We employ data from Italian Survey of Household Income and Wealth (SHIW), collected by the Bank of Italy \cite{bankofitaly09}. For the purpose of our study, from the sampled records of about 8,000 households annual income distributions in all 20 Italian regions were aggregated into form of distribution-valued variables. In order to exclude outlying observations that would destroy the smoothing process, from non-zero income values in the whole data set those above 99\% quantile were omitted. Consequently, for values from single regions the optimal number of classes was computed (according to the well-known Sturges rule), resulting in mean value $9.24$, i.e. nine equidistant income classes in the range $I=[65.7,1.10709\times 10^5]$ were constructed. Data $x_i,\,i=1,\dots,9$ thus correspond to midpoints of income intervals, relative frequencies (proportions of income classes on the overall income distribution in single regions) stand for $y_i$ values. Because of count character of values in income classes, also some zero values occurred that would make further processing using centred logratio transformation not possible. For this reason, their imputation using a model-based procedure was performed \cite{martin14}; the resulting complete data are collected in Table \ref{propor}.
Moreover, for subsequent interpretation purposes, the employed Italian regions were divided into three parts according to their geographical location ($\mbox{north}=\mbox{N}$, $\mbox{middle}=\mbox{M}$, $\mbox{south and islands}=\mbox{S}$), following the National Statistical Institute (ISTAT) classification.

\setcounter{table}{0}

\begin{table}
\tbl{Proportions of income classes in 20 Italian regions.}
{\begin{tabular}[l]{@{}lc ccccccccc }\toprule
    \textit{region} & \textit{loc.} & \multicolumn{9}{c}{\textit{proportions of income classes, $\quad y_i$, $ i=1,\dots,9$}} \\
    \toprule
    Piemonte & N & 0.067 & 0.385 & 0.323 & 0.134 & 0.052 & 0.022 & 0.009 & 0.005 & 0.003 \\
    Valle d'Aosta & N & 0.042 & 0.340 & 0.319 & 0.212 & 0.042 & 0.016 & 0.006 & 0.016 & 0.006 \\
    Lombardia & N & 0.089 & 0.275 & 0.269 & 0.151 & 0.107 & 0.056 & 0.022 & 0.018 & 0.012 \\
    Trentino & N & 0.058 & 0.320 & 0.279 & 0.127 & 0.134 & 0.029 & 0.041 & 0.006 & 0.005 \\
    Veneto & N & 0.103 & 0.329 & 0.255 & 0.177 & 0.081 & 0.022 & 0.015 & 0.010 & 0.007 \\
    Friuli & N & 0.084 & 0.320 & 0.232 & 0.168 & 0.088 & 0.068 & 0.028 & 0.008 & 0.004 \\
    Liguria & N & 0.081 & 0.362 & 0.207 & 0.213 & 0.081 & 0.026 & 0.026 & 0.003 & 0.002 \\
    Emilia Romagna & N & 0.065 & 0.303 & 0.275 & 0.189 & 0.085 & 0.045 & 0.017 & 0.015 & 0.006 \\
    Toscana & M & 0.043 & 0.283 & 0.293 & 0.188 & 0.105 & 0.052 & 0.015 & 0.015 & 0.007 \\
    Umbria & M & 0.052 & 0.351 & 0.337 & 0.157 & 0.056 & 0.026 & 0.015 & 0.004 & 0.002 \\
    Marche & M & 0.115 & 0.401 & 0.219 & 0.153 & 0.058 & 0.032 & 0.014 & 0.006 & 0.003 \\
    Lazio & M & 0.120 & 0.349 & 0.260 & 0.150 & 0.066 & 0.032 & 0.012 & 0.007 & 0.002 \\
    Abruzzo & S & 0.100 & 0.368 & 0.294 & 0.144 & 0.045 & 0.030 & 0.004 & 0.010 & 0.005 \\
    Molise & S & 0.131 & 0.349 & 0.277 & 0.109 & 0.080 & 0.022 & 0.022 & 0.006 & 0.004 \\
    Campania & S & 0.238 & 0.485 & 0.167 & 0.066 & 0.019 & 0.016 & 0.006 & 0.002 & 0.001 \\
    Puglia & S & 0.238 & 0.441 & 0.201 & 0.068 & 0.025 & 0.009 & 0.011 & 0.003 & 0.002 \\
    Basilicata & S & 0.246 & 0.385 & 0.169 & 0.115 & 0.038 & 0.031 & 0.006 & 0.006 & 0.003 \\
    Calabria & S & 0.240 & 0.408 & 0.209 & 0.084 & 0.037 & 0.005 & 0.010 & 0.004 & 0.003 \\
    Sicilia & S & 0.255 & 0.473 & 0.161 & 0.053 & 0.029 & 0.014 & 0.012 & 0.002 & 0.001 \\
    Sardegna & S & 0.167 & 0.425 & 0.217 & 0.123 & 0.044 & 0.015 & 0.006 & 0.003 & 0.002 \\
    \toprule
    interval midpoints & & 6574 & 19591 & 32608 & 45625 & 58641 & 71658 & 84675 & 97692 & 110709 \\
    \botrule
\end{tabular}}
\label{propor}
\end{table}

Furthermore, we perform the discrete version of clr transformation \cite{aitchison86}, defined as $z_i=\ln\frac{y_i}{g(y_1,\dots,y_9)}$, where $g$ denotes geometric mean of the argument values; as expected, a condition $\sum_{i=1}^9z_i=0$ holds. The obtained values for all employed regions are displayed in Table \ref{clrincome}.

\begin{table}
\tbl{Centred logratio transformation of income classes' proportions in 20 Italian regions.}
{\begin{tabular}[l]{@{}lc rrrrrrrrr }\toprule
    \textit{region} & \textit{loc.} & \multicolumn{9}{c}{\textit{clr transformation of income classes' proportions, $\quad z_i$, $i=1,\dots,9$}}\\
    \toprule
    Piemonte & N & 0.587 & 2.331 & 2.154 & 1.271 & 0.331 & -0.550 & -1.437 & -1.997 & -2.690 \\
    Valle d'Aosta & N & 0.015 & 2.094 & 2.030 & 1.624 & 0.015 & -0.946 & -1.919 & -0.966 & -1.946 \\
    Lombardia & N & 0.275 & 1.405 & 1.383 & 0.805 & 0.462 & -0.187 & -1.125 & -1.307 & -1.713 \\
    Trentino & N & 0.084 & 1.789 & 1.653 & 0.873 & 0.917 & -0.609 & -0.272 & -2.218 & -2.218 \\
    Veneto & N & 0.681 & 1.843 & 1.588 & 1.224 & 0.442 & -0.865 & -1.232 & -1.638 & -2.043 \\
    Friuli & N & 0.401 & 1.739 & 1.417 & 1.095 & 0.448 & 0.190 & -0.697 & -1.950 & -2.643 \\
    Liguria & N & 0.666 & 2.166 & 1.606 & 1.637 & 0.666 & -0.473 & -0.473 & -2.662 & -3.132 \\
    Emilia Romagna & N & 0.140 & 1.681 & 1.584 & 1.209 & 0.405 & -0.223 & -1.204 & -1.291 & -2.303 \\
    Toscana & M & -0.259 & 1.619 & 1.654 & 1.211 & 0.627 & -0.083 & -1.319 & -1.319 & -2.130 \\
    Umbria & M & 0.348 & 2.252 & 2.209 & 1.447 & 0.417 & -0.345 & -0.905 & -2.291 & -3.133 \\
    Marche & M & 0.965 & 2.211 & 1.607 & 1.247 & 0.272 & -0.326 & -1.114 & -2.031 & -2.831 \\
    Lazio & M & 0.982 & 2.046 & 1.753 & 1.201 & 0.386 & -0.345 & -1.301 & -1.811 & -2.910 \\
    Abruzzo & S & 0.856 & 2.164 & 1.938 & 1.227 & 0.057 & -0.348 & -2.308 & -1.447 & -2.140 \\
    Molise & S & 1.000 & 1.982 & 1.748 & 0.818 & 0.508 & -0.791 & -0.791 & -2.071 & -2.404 \\
    Campania & S & 2.267 & 2.980 & 1.914 & 0.983 & -0.245 & -0.428 & -1.344 & -2.730 & -3.396 \\
    Puglia & S & 2.009 & 2.628 & 1.844 & 0.756 & -0.247 & -1.258 & -1.035 & -1.952 & -2.746 \\
    Basilicata & S & 1.839 & 2.286 & 1.465 & 1.082 & -0.017 & -0.240 & -1.841 & -1.933 & -2.641 \\
    Calabria & S & 1.988 & 2.516 & 1.848 & 0.931 & 0.105 & -1.841 & -1.148 & -2.100 & -2.298 \\
    Sicilia & S & 2.145 & 2.763 & 1.684 & 0.574 & -0.014 & -0.776 & -0.931 & -2.722 & -2.722 \\
    Sardegna & S & 1.657 & 2.591 & 1.918 & 1.351 & 0.322 & -0.777 & -1.693 & -2.521 & -2.848 \\
    \toprule
    interval midpoints & & 6574 & 19591 & 32608 & 45625 & 58641 & 71658 & 84675 & 97692 & 110709 \\
    \botrule
\end{tabular}}
\label{clrincome}
\end{table}

Now we can finally proceed to approximate the clr transformed proportions in the $L^2$ space with optimal smoothing splines from Section 3. Let $\mathbf{x}$ denotes the vector of interval midpoints and $w_{i}=1$ for $i=1,\dots,9$. For the purpose of our study, the parameter $\alpha$ was set simply to $1$, other options are discussed in \cite{deboor78}. 
Furthermore, we set $k=3$, $l=2$ (i.e., cubic splines are employed) and we will consider knots
$\Delta\lambda:=0<30\,000<70\,000<1.10709\times 10^5$ for our optimal smoothing problem. The corresponding spline coefficients are
collected in Table \ref{clrcoef}; note that the zero sum constraint is fulfilled as well.

\begin{figure}[h]
\begin{center}
\includegraphics[width=14cm]{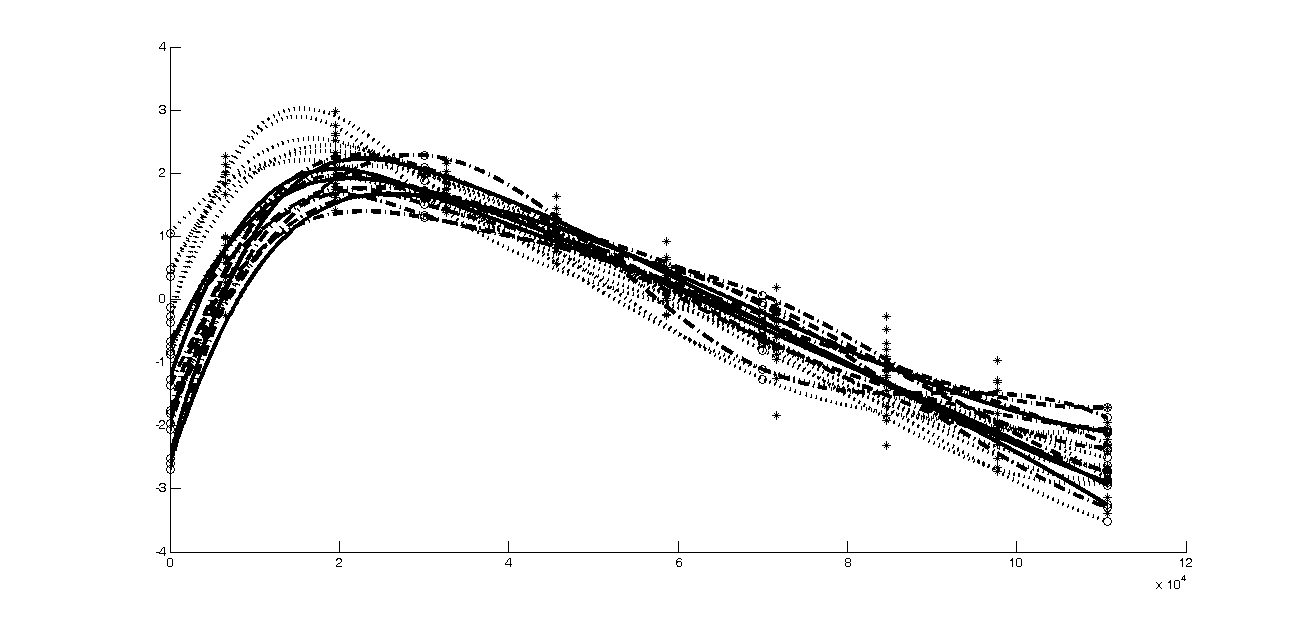}
\includegraphics[width=14cm]{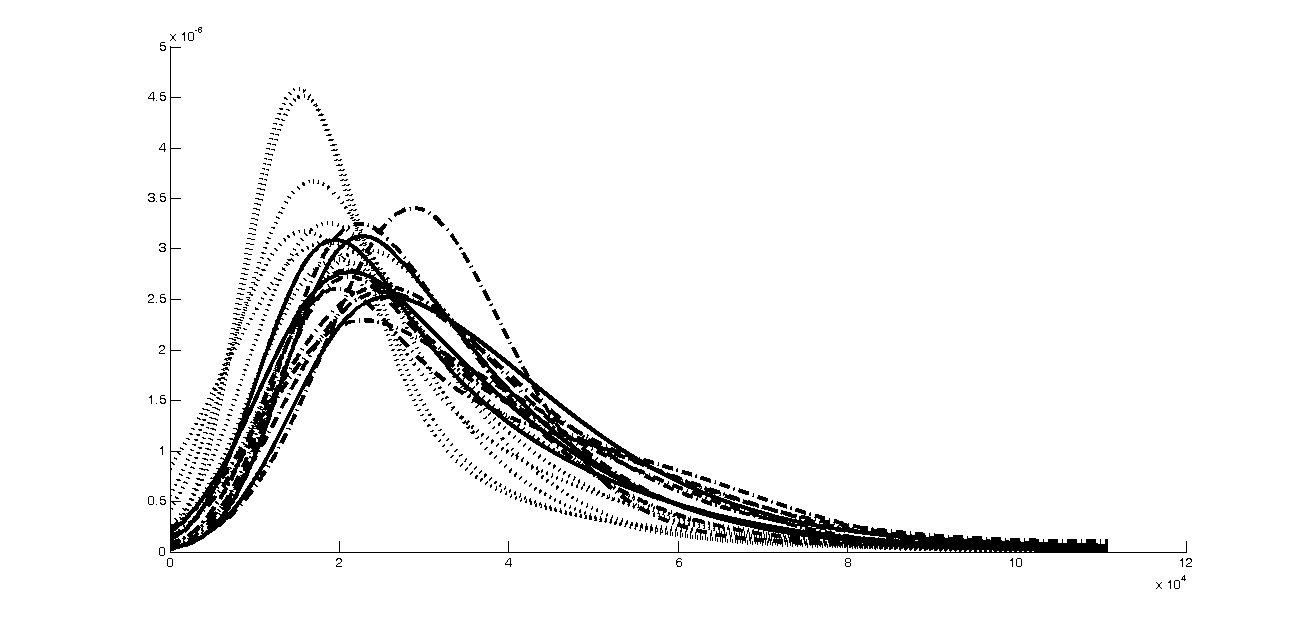}
\caption{Optimal smoothing spline for SHIW data for all regions for clr transformed data (up) and after its inverse clr transformation (down). Lines are denoted according to geographical location of regions: N - dashed black line, M - solid black line, S - dotted grey line.}\label{f1}
\end{center}
\end{figure}

The resulting optimal smoothing splines can be also back-transformed to the original Bayes space of density functions, see Figure \ref{f1}. It is easy to see from the up plot that the $B$-spline functions nicely follow the clr-transformed proportions. This is due to both the proper underlying geometry for the approximation procedure, provided by the $L^2$ space, but also the choice of knots that avoids nuisance variability of densities. Moreover, regional effects of income distribution are nice visible, in particular different distributions for northern and middle Italian regions comparing to the southern ones.

\begin{table}
\tbl{B-spline coefficients for clr transformed density functions of 20 Italian regions.}
{\begin{tabular}[l]{@{}lc rrrrrr }
    \toprule
    \textit{region} & \textit{loc.} & \multicolumn{6}{c}{\textit{spline coefficients for given knots, $b_i^*,\, i=-3,\dots,2$}} \\
    \toprule
    Piemonte & N & 1.972 & 2.650 & 2.376 & -0.907 & -2.202 & -2.734 \\
    Valle d'Aosta & N & -1.801 & 1.192 & 3.979 & -2.791 & -1.048 & -1.877 \\
    Lombardia & N & -1.362 & 1.609 & 1.413 & -0.129 & -1.788 & -1.708 \\
    Trentino & N & -2.686 & 2.512 & 1.128 & 0.519 & -2.250 & -2.357 \\
    Veneto & N & -0.660 & 1.602 & 2.370 & -1.079 & -1.872 & -2.067 \\
    Friuli & N & -2.065 & 2.531 & 0.742 & 0.828 & -2.205 & -2.728 \\
    Liguria & N & -1.756 & 2.628 & 1.467 & 0.487 & -2.643 & -3.299 \\
    Emilia Romagna & N & -1.806 & 1.609 & 2.080 & -0.666 & -1.506 & -2.297 \\
    Toscana & M & -2.595 & 1.615 & 2.038 & -0.254 & -1.823 & -2.101 \\
    Umbria & M & -2.517 & 2.625 & 2.226 & -0.444 & -2.144 & -3.254 \\
    Marche & M & -1.260 & 2.816 & 1.418 & -0.160 & -2.239 & -2.896 \\
    Lazio & M & -0.750 & 2.247 & 1.889 & -0.555 & -2.015 & -2.944 \\
    Abruzzo & S & 0.872 & 2.184 & 2.542 & -1.268 & -2.280 & -2.056 \\
    Molise & S & -0.827 & 2.474 & 1.541 & -0.448 & -2.073 & -2.508 \\
    Campania & S & -0.275 & 4.574 & 0.753 & -0.177 & -2.878 & -3.523 \\
    Puglia & S & 0.372 & 3.304 & 1.843 & -1.926 & -1.449 & -2.856 \\
    Basilicata & S & 0.476 & 2.974 & 1.174 & -0.267 & -2.672 & -2.633 \\
    Calabria & S & 1.041 & 2.570 & 2.561 & -2.272 & -1.802 & -2.400 \\
    Sicilia & S & -0.360 & 4.563 & 0.260 & 0.097 & -2.864 & -2.874 \\
    Sardegna & S & -0.124 & 3.085 & 1.937 & -0.529 & -3.097 & -2.903\\
\end{tabular}}
\label{clrcoef}
\end{table}

\begin{figure}[h]
\centering
\includegraphics[width=14cm]{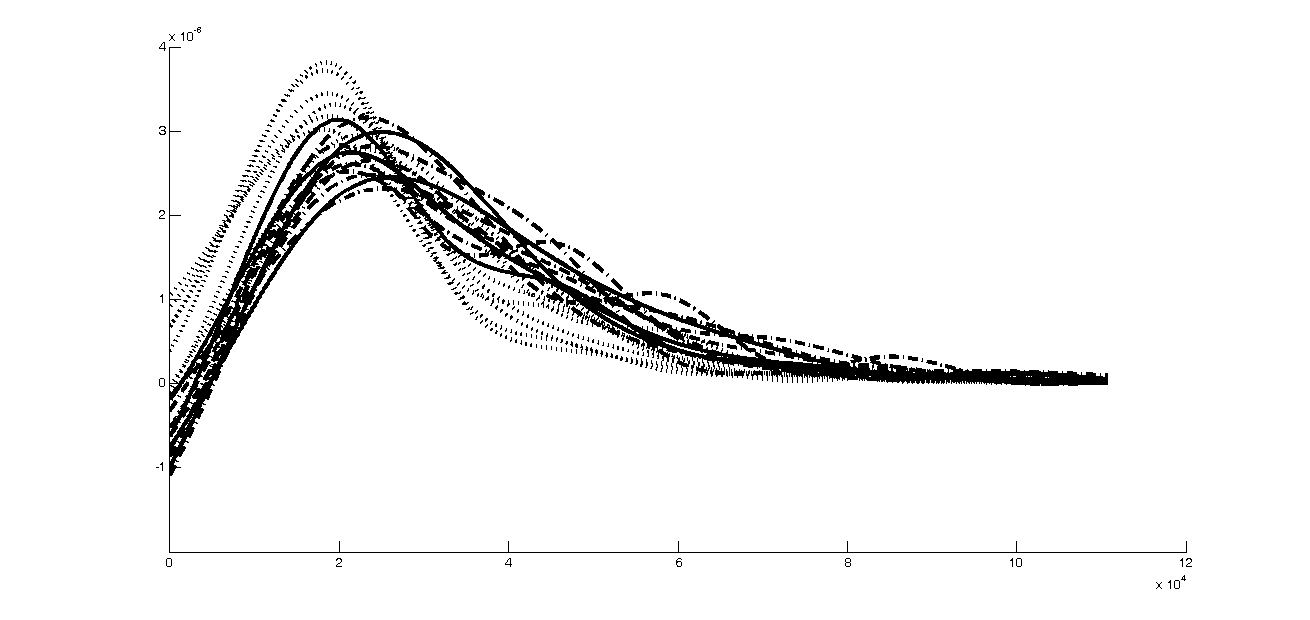}
\caption{Cubic smoothing spline for SHIW data for all regions using the original proportional data.}\label{f2}
\end{figure}

\begin{figure}[h]
\centering
\includegraphics[width=14cm]{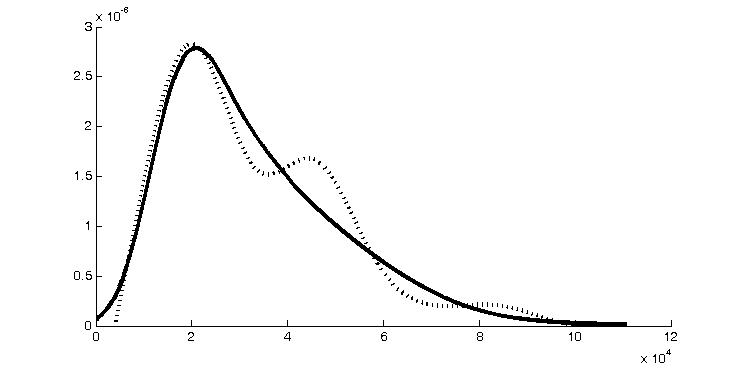}
\includegraphics[width=14cm]{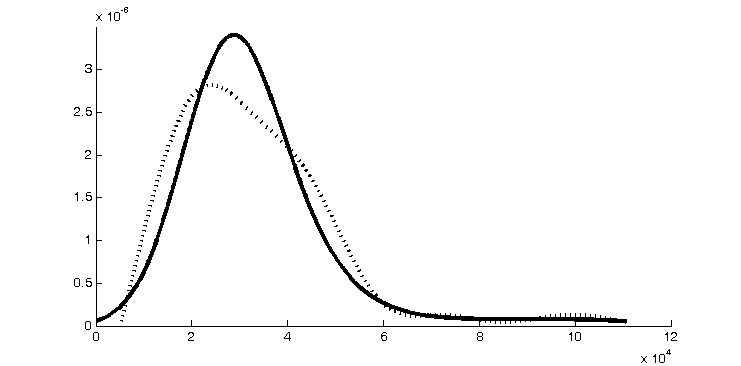}
\caption{Comparison between the presented approach (solid black line) and standard smoothing problem (dotted grey line) for Liguria (up) and Valle d'Aosta (down); only positive values of splines are displayed.}\label{f3}
\end{figure}

Finally, we compare the presented approach with an alternative concept of using cubic smoothing splines for the original proportional data (Table \ref{propor}) as proposed in \cite{delicado11}, where the knots correspond to data points $x_i$, $i=1,\dots,9$; furthermore, parameters $k=3$ and $l=2$ were set as before. From Figure \ref{f2} with the resulting smoothed densities several features can be easily derived. In addition to choppy right tails as a consequence of inappropriate choice of knots (a typical situation is shown in Figure \ref{f3}, up, for the case of Liguria region) also negative values are reached on the left tail of curves, and for some of regions even on the right tail of densities as well. Although numerically this undesired effect (that was automatically avoided before) could be overcome by the logarithmic transformation of the original density functions and subsequent maximum-likelihood estimation of the $B$-spline coefficients as proposed in \cite{ramsay05}, this approach would still ignore the inherent geometrical features of densities, captured by the Bayes space methodology. Moreover, as a consequence of the relative scale property of density functions, captured by smoothing in the clr space, also data structure itself seems to be quite different for both approaches. For example, Valle d'Aosta with the highest maximum among income distributions of northern regions (see Figure \ref{f1}, down) shows an exceptional behaviour, when the estimation is performed in the clr space. This could be easily explained by the fact that in this region there is indeed a lot of wealth with respect to the rest of Italy and the local population gains also from a low-tax policy due to the special autonomous status of the Valle d'Aosta region. On the other hand, this particular behaviour is rather not visible, when approximating the original proportions.

\section{Conclusions}


With modern monitoring tools that enable to produce huge amounts of data, statistical analysis of density functions become more and more important. Although the theoretical background for their meaningful statistical treatment, provided by Bayes spaces, is already available, for its practical applicability a wide spectrum of mathematical problems needs to be solved. This paper handled one of them, smoothing of clr transformed density functions. We have shown that $B$-splines represent an easy-to-handle tool for this purpose, moreover, the special case of smoothing splines provides a trade-off between the interpolation problem and the least squares approximation. For the real-world data example, where income distributions in 20 Italian regions were analyzed, the advantages of the Bayes space approach were clearly demonstrated by comparison with approximation of the original densities using cubic smoothing splines as proposed in \cite{delicado11}. In addition to the above advances also ``dimention reduction'' of the discrete input data, provided by the $B$-spline coefficients, should be mentioned that enables to analyze statistically also high-dimensional functional data (resulting, e.g., from metabolomical experiments). We avoided to compare the results with the case of approximation of density functions with Bernain polynomials, proposed in \cite{menafoglio14}, as they represent a completely different concept of representation of densities; nevertheless, also here similar effects, resulting from ignoring the relative character of density functions, could be observed.

Of course, this paper cannot be considered as a final solution, it rather opens a further variety of challenges concerning smoothing of density functions using Bayes spaces (concerning choice of basis functions, optimal choice of knots etc.). On the other hand, we are convinced that it provides a concise approach and set a clear direction for future research developments in the field.



\label{lastpage}

\def\refname{References}


%

\end{document}